\documentclass[12pt]{article}%
\usepackage{amsfonts}
\usepackage{sw20bams}%
\usepackage{amsmath}%
\usepackage{amssymb}%
\usepackage{graphicx}
\providecommand{\U}[1]{\protect\rule{.1in}{.1in}}
\begin{document}

\title{Union of $n$ Disks:\ Remote Centers, Common Origin}
\author{Steven Finch}
\date{November 16, 2015}
\maketitle

\begin{abstract}
Explicit area expressions are known for a special case, due to Tao \&\ Wu
(1987), and lead to calculation of integrals in applied probability.

\end{abstract}

\footnotetext{Copyright \copyright \ 2015 by Steven R. Finch. All rights
reserved.}A collection of $n$ planar disks is said to have \textbf{remote
centers} if the $i^{\text{th}}$ disk never contains the $j^{\text{th}}$center,
for any $1\leq i\neq j\leq n$. \ It has \textbf{common origin} if the point
$\vec{0}$ is on the boundary of each disk. \ No further assumptions are made
concerning the relative sizes of the disks or the extent to which they
overlap. \ Let the disks be centered at $\vec{r}_{1}$, $\vec{r}_{2}$, \ldots,
$\vec{r}_{n}$. \ It follows that $\vec{0}$ is closer to $\vec{r}_{j}$ than any
other $\vec{r}_{i}$. Let $r_{j}$ denote the length of $\vec{r}_{j}$ (the
$j^{\text{th}}$ radius) and $r_{i,j}$ denote the length of $\vec{r}_{i}-$
$\vec{r}_{j}$. \ The constraints $r_{j}<r_{i,j}$ are crucial for the
calculation of certain integrals in \cite{TW-DsUn}, which in turn give the
probability that an individual survives a random type of violent shootout.

Consider the integral%
\[
c_{n}=\frac{1}{n!}%
{\displaystyle\int\limits_{r_{j}<r_{i,j}}}
\exp\left[  -V\left(  \vec{r}_{1},\vec{r}_{2},\ldots,\vec{r}_{n}\right)
\right]  d\vec{r}_{1}d\vec{r}_{2}\ldots d\vec{r}_{n}
\]
where $V\left(  \vec{r}_{1},\vec{r}_{2},\ldots,\vec{r}_{n}\right)  $ is the
area of the union of $n$ disks centered at $\vec{r}_{1}$, $\vec{r}_{2}$,
\ldots, $\vec{r}_{n}$ and intersecting at $\vec{0}$. \ The value of $V$ can be
found symbolically by use of RegionUnion and RegionMeasure functions in
\textsc{Mathematica 10}. \ Such technology does not currently allow us to
evaluate $c_{n}$ directly because computer memory is quickly exhausted in the
required numerical quadrature. \ One aim of this paper is to revisit formulas
in \cite{TW-DsUn}, focusing on an elaborate change of variables and leading to
a less-intensive indirect calculation. \ The notation employed previously is
adopted here too, so that several details can be clarified.

Another aim of this paper is to give explicit expressions for $V$. \ We are
aware of a substantial literature devoted to the problem of arbitrary
unions/intersections of disks/balls \cite{X1-DsUn, X2-DsUn, X3-DsUn, X4-DsUn,
X5-DsUn, X6-DsUn, X7-DsUn, X8-DsUn, X9-DsUn, X10-DsUn, X11-DsUn, X12-DsUn,
X13-DsUn, X14-DsUn, X15-DsUn}, yet have not seen anything (in an unsystematic
survey) resembling Tao \&\ Wu's results. \ It is important to remember that
our disks \textit{always} possess remote centers and common origin. \ Thus our
results are not general, but nevertheless might constitute an interesting
special case for future study.

\section{Areas}

We assume WLOG that the $n$ disk centers are sorted according to increasing
argument (counterclockwise angle with respect to the horizontal axis). Since
$V$ is homogeneous and quadratic in $\vec{r}_{1}$, $\vec{r}_{2}$, \ldots,
$\vec{r}_{n}$, we can reparametrize it as follows:
\[
V=r_{1}^{2}V\left(  t_{2},t_{3},\ldots,t_{n};\theta_{1},\theta_{2}%
,\ldots,\theta_{n-1}\right)
\]
where $t_{k}=r_{k}/r_{k-1}$ and $\theta_{k}$ is the (positive) angle between
$\vec{r}_{k}$ and $\vec{r}_{k+1}$. \ Let $\theta_{n}=2\pi-\theta_{1}%
-\theta_{2}-\cdots-\theta_{n-1}$, which is well-defined (positive) by the
ordering in our construction. \ For $n=2$,
\[
V\left(  t_{2};\theta_{1}\right)  =z_{1}+t_{2}^{2}z_{2}%
\]
where%
\[
z_{1}=\pi-\alpha_{1}+\frac{1}{2}\sin(2\alpha_{1})-\beta_{2}+\frac{1}{2}%
\sin(2\beta_{2}),
\]%
\[
z_{2}=\pi-\alpha_{2}+\frac{1}{2}\sin(2\alpha_{2})-\beta_{1}+\frac{1}{2}%
\sin(2\beta_{1})
\]
and%
\[
\sin(\alpha_{1})=\left\{
\begin{array}
[c]{lll}%
\dfrac{t_{2}\sin(\theta_{1})}{\sqrt{1+t_{2}^{2}-2t_{2}\cos(\theta_{1})}} &  &
\text{if }\theta_{1}<\pi,\\
0 &  & \text{if }\theta_{1}>\pi\text{;}%
\end{array}
\right.
\]%
\[
\sin(\beta_{1})=\left\{
\begin{array}
[c]{lll}%
\dfrac{\sin(\theta_{1})}{\sqrt{1+t_{2}^{2}-2t_{2}\cos(\theta_{1})}} &  &
\text{if }\theta_{1}<\pi,\\
0 &  & \text{if }\theta_{1}>\pi\text{;}%
\end{array}
\right.
\]%
\[
\sin(\alpha_{2})=\left\{
\begin{array}
[c]{lll}%
\dfrac{\sin(\theta_{2})}{\sqrt{1+t_{2}^{2}-2t_{2}\cos(\theta_{2})}} &  &
\text{if }\theta_{2}<\pi,\\
0 &  & \text{if }\theta_{2}>\pi\text{;}%
\end{array}
\right.
\]%
\[
\sin(\beta_{2})=\left\{
\begin{array}
[c]{lll}%
\dfrac{t_{2}\sin(\theta_{2})}{\sqrt{1+t_{2}^{2}-2t_{2}\cos(\theta_{2})}} &  &
\text{if }\theta_{2}<\pi,\\
0 &  & \text{if }\theta_{2}>\pi\text{.}%
\end{array}
\right.
\]
The preceding is needlessly complicated. \ Assume additionally that
$\theta_{1}<\pi$ as pictured in Figure 1; it follows that $\theta_{2}%
=2\pi-\theta_{1}>\pi$, $\alpha_{2}=0$, \ $\beta_{2}=0$ and therefore%
\[%
\begin{array}
[c]{ccc}%
z_{1}=\pi-\alpha_{1}+\dfrac{1}{2}\sin(2\alpha_{1}), &  & z_{2}=\pi-\beta
_{1}+\dfrac{1}{2}\sin(2\beta_{1}).
\end{array}
\]
\{Correction to \cite{TW-DsUn}: angle $\beta_{2}$ should be $\beta_{1}$ in
formula (22) and angle $\beta_{i+1}$ should be $\beta_{i}$ in formula (23).
\ The corresponding figure, however, is fine.\}%
\begin{figure}[ptb]%
\centering
\includegraphics[
height=2.8997in,
width=3.1176in
]%
{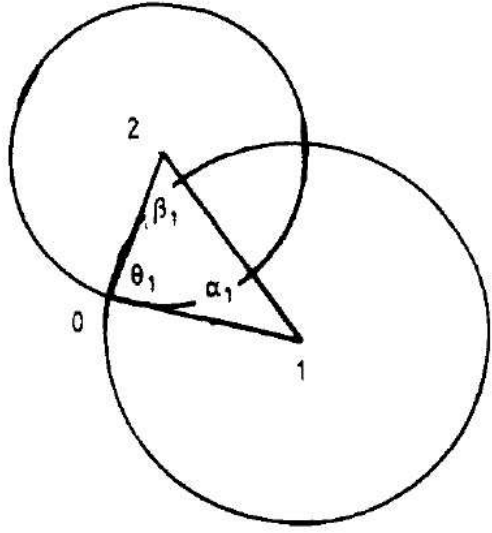}%
\caption{Area occupied by two circles (from \cite{TW-DsUn}).}%
\end{figure}

For $n=3$,
\[
V\left(  t_{2},t_{3};\theta_{1},\theta_{2}\right)  =z_{1}+t_{2}^{2}z_{2}%
+t_{2}^{2}t_{3}^{2}z_{3}%
\]
where%
\[
z_{1}=\pi-\alpha_{1}+\frac{1}{2}\sin(2\alpha_{1})-\beta_{3}+\frac{1}{2}%
\sin(2\beta_{3}),
\]%
\[
z_{2}=\pi-\alpha_{2}+\frac{1}{2}\sin(2\alpha_{2})-\beta_{1}+\frac{1}{2}%
\sin(2\beta_{1}),
\]%
\[
z_{3}=\pi-\alpha_{3}+\frac{1}{2}\sin(2\alpha_{3})-\beta_{2}+\frac{1}{2}%
\sin(2\beta_{2})
\]
and%
\[
\sin(\alpha_{1})=\left\{
\begin{array}
[c]{lll}%
\dfrac{t_{2}\sin(\theta_{1})}{\sqrt{1+t_{2}^{2}-2t_{2}\cos(\theta_{1})}} &  &
\text{if }\theta_{1}<\pi,\\
0 &  & \text{if }\theta_{1}>\pi\text{;}%
\end{array}
\right.
\]%
\[
\sin(\beta_{1})=\left\{
\begin{array}
[c]{lll}%
\dfrac{\sin(\theta_{1})}{\sqrt{1+t_{2}^{2}-2t_{2}\cos(\theta_{1})}} &  &
\text{if }\theta_{1}<\pi,\\
0 &  & \text{if }\theta_{1}>\pi\text{;}%
\end{array}
\right.
\]%
\[
\sin(\alpha_{2})=\left\{
\begin{array}
[c]{lll}%
\dfrac{t_{3}\sin(\theta_{2})}{\sqrt{1+t_{3}^{2}-2t_{3}\cos(\theta_{2})}} &  &
\text{if }\theta_{2}<\pi,\\
0 &  & \text{if }\theta_{2}>\pi\text{;}%
\end{array}
\right.
\]%
\[
\sin(\beta_{2})=\left\{
\begin{array}
[c]{lll}%
\dfrac{\sin(\theta_{2})}{\sqrt{1+t_{3}^{2}-2t_{3}\cos(\theta_{2})}} &  &
\text{if }\theta_{2}<\pi,\\
0 &  & \text{if }\theta_{2}>\pi\text{;}%
\end{array}
\right.
\]%
\[
\sin(\alpha_{3})=\left\{
\begin{array}
[c]{lll}%
\dfrac{\sin(\theta_{3})}{\sqrt{1+t_{2}^{2}t_{3}^{2}-2t_{2}t_{3}\cos(\theta
_{3})}} &  & \text{if }\theta_{3}<\pi,\\
0 &  & \text{if }\theta_{3}>\pi\text{;}%
\end{array}
\right.
\]%
\[
\sin(\beta_{3})=\left\{
\begin{array}
[c]{lll}%
\dfrac{t_{2}t_{3}\sin(\theta_{3})}{\sqrt{1+t_{2}^{2}t_{3}^{2}-2t_{2}t_{3}%
\cos(\theta_{3})}} &  & \text{if }\theta_{3}<\pi,\\
0 &  & \text{if }\theta_{3}>\pi\text{.}%
\end{array}
\right.
\]
Figure 2 provides a sample scenario in which the overlap is merely the point
$\vec{0}$; a more representative picture would include a nondegenerate
(circular triangle) intersection. \{The conditions for $\alpha_{k}$ and
$\beta_{k}$ to vanish are not fully stated in \cite{TW-DsUn}; we have
attempted to be more precise here.\}%
\begin{figure}[ptb]%
\centering
\includegraphics[
height=4.3405in,
width=3.8744in
]%
{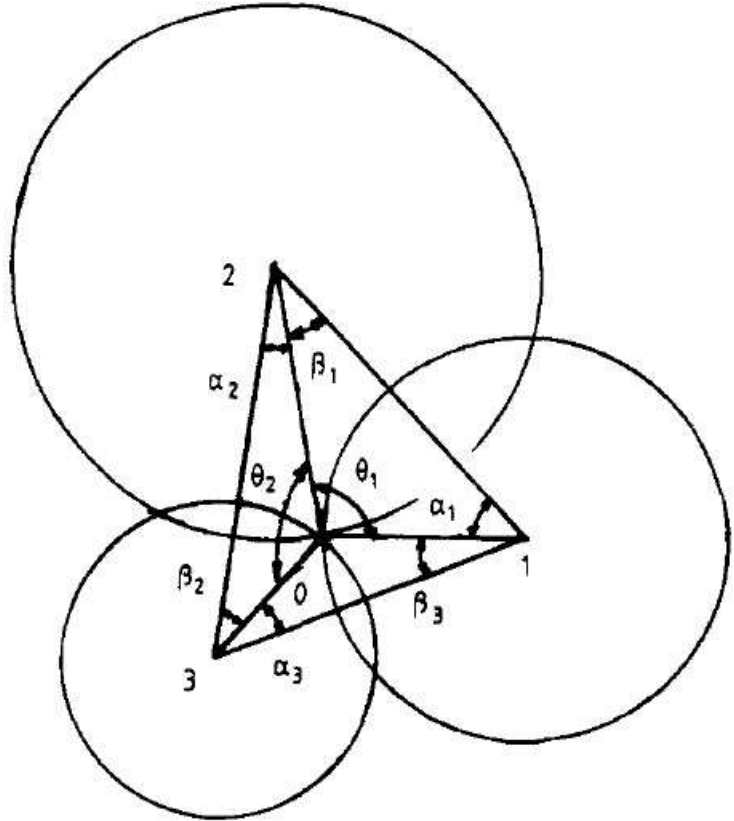}%
\caption{Area occupied by three circles intersecting at a single point (from
\cite{TW-DsUn}).}%
\end{figure}

For arbitrary $n$,
\[
V\left(  t_{2},t_{3},\ldots,t_{n};\theta_{1},\theta_{2},\ldots,\theta
_{n-1}\right)  =z_{1}+t_{2}^{2}z_{2}+t_{2}^{2}t_{3}^{2}z_{3}+\cdots+t_{2}%
^{2}t_{3}^{2}\cdots t_{n}^{2}z_{n}
\]
where%
\[
z_{k}=\pi-\alpha_{k}+\frac{1}{2}\sin(2\alpha_{k})-\beta_{k-1}+\frac{1}{2}%
\sin(2\beta_{k-1})
\]
for $1\leq k\leq n$ and we agree to set $\beta_{0}=\beta_{n}$. \ Also,%
\[
\sin(\alpha_{k})=\left\{
\begin{array}
[c]{lll}%
\dfrac{t_{k+1}\sin(\theta_{k})}{\sqrt{1+t_{k+1}^{2}-2t_{k+1}\cos(\theta_{k})}}
&  & \text{if }\theta_{k}<\pi,\\
0 &  & \text{if }\theta_{k}>\pi\text{;}%
\end{array}
\right.
\]%
\[
\sin(\beta_{k})=\left\{
\begin{array}
[c]{lll}%
\dfrac{\sin(\theta_{k})}{\sqrt{1+t_{k+1}^{2}-2t_{k+1}\cos(\theta_{k})}} &  &
\text{if }\theta_{k}<\pi,\\
0 &  & \text{if }\theta_{k}>\pi
\end{array}
\right.
\]
for $1\leq k\leq n-1$ and
\[
\sin(\alpha_{n})=\left\{
\begin{array}
[c]{lll}%
\dfrac{\sin(\theta_{n})}{\sqrt{1+t_{2}^{2}t_{3}^{2}\cdots t_{n}^{2}%
-2t_{2}t_{3}\cdots t_{n}\cos(\theta_{n})}} &  & \text{if }\theta_{n}<\pi,\\
0 &  & \text{if }\theta_{n}>\pi\text{;}%
\end{array}
\right.
\]%
\[
\sin(\beta_{n})=\left\{
\begin{array}
[c]{lll}%
\dfrac{t_{2}t_{3}\cdots t_{n}\sin(\theta_{n})}{\sqrt{1+t_{2}^{2}t_{3}%
^{2}\cdots t_{n}^{2}-2t_{2}t_{3}\cdots t_{n}\cos(\theta_{n})}} &  & \text{if
}\theta_{n}<\pi,\\
0 &  & \text{if }\theta_{n}>\pi\text{.}%
\end{array}
\right.
\]
A proof of these formulas does not appear in \cite{TW-DsUn}. \ We point out
only that, for $1\leq k\leq n-1$,
\[
\frac{\sin(\alpha_{k})}{r_{k+1}}=\frac{\sin(\theta_{k})}{r_{k,k+1}}=\frac
{\sin(\beta_{k})}{r_{k}}
\]
by the Law of Sines,%
\[
r_{k,k+1}^{2}=r_{k}^{2}+r_{k+1}^{2}-2r_{k}r_{k+1}\cos(\theta_{k})
\]
by the Law of Cosines, $t_{k+1}=r_{k+1}/r_{k}$, hence%
\[
\frac{\sin(\alpha_{k})}{t_{k+1}}=\frac{\sin(\alpha_{k})}{r_{k+1}/r_{k}}%
=\frac{\sin(\theta_{k})}{r_{k,k+1}/r_{k}}=\sin(\beta_{k}).
\]
Also%
\[
\frac{\sin(\alpha_{n})}{r_{1}}=\frac{\sin(\theta_{n})}{r_{1,n}}=\frac
{\sin(\beta_{n})}{r_{n}}
\]
by the Law of Sines,%
\[
r_{1,n}^{2}=r_{1}^{2}+r_{n}^{2}-2r_{1}r_{n}\cos(\theta_{n})
\]
by the Law of Cosines, $t_{2}t_{3}\cdots t_{n}=r_{n}/r_{1}$, hence%
\[
\sin(\alpha_{n})=\frac{\sin(\theta_{n})}{r_{1,n}/r_{1}}=\frac{\sin(\beta_{n}%
)}{r_{n}/r_{1}}=\frac{\sin(\beta_{n})}{t_{2}t_{3}\cdots t_{n}}.
\]

\section{Integrals}

Let $\vec{r}_{j}=(x_{j},y_{j})$ for each $1\leq j\leq n$ and assume that the
argument of $\vec{r}_{1}$ is $\theta_{0}$.\ The change of variables%
\begin{align*}
\left(
\begin{array}
[c]{c}%
x_{1}\\
y_{1}\\
x_{2}\\
y_{2}\\
x_{3}\\
y_{3}\\
\vdots\\
x_{n}\\
y_{n}%
\end{array}
\right)   & =\left(
\begin{array}
[c]{c}%
r_{1}\cos(\theta_{0})\\
r_{1}\sin(\theta_{0})\\
r_{2}\cos(\theta_{0}+\theta_{1})\\
r_{2}\sin(\theta_{0}+\theta_{1})\\
r_{3}\cos(\theta_{0}+\theta_{1}+\theta_{2})\\
r_{3}\sin(\theta_{0}+\theta_{1}+\theta_{2})\\
\vdots\\
r_{n}\cos(\theta_{0}+\theta_{1}+\theta_{2}+\cdots+\theta_{n-1})\\
r_{n}\sin(\theta_{0}+\theta_{1}+\theta_{2}+\cdots+\theta_{n-1})
\end{array}
\right) \\
& =r_{1}\left(
\begin{array}
[c]{c}%
\cos(\theta_{0})\\
\sin(\theta_{0})\\
t_{2}\cos(\theta_{0}+\theta_{1})\\
t_{2}\sin(\theta_{0}+\theta_{1})\\
t_{2}t_{3}\cos(\theta_{0}+\theta_{1}+\theta_{2})\\
t_{2}t_{3}\sin(\theta_{0}+\theta_{1}+\theta_{2})\\
\vdots\\
t_{2}t_{3}\cdots t_{n}\cos(\theta_{0}+\theta_{1}+\theta_{2}+\cdots
+\theta_{n-1})\\
t_{2}t_{3}\cdots t_{n}\sin(\theta_{0}+\theta_{1}+\theta_{2}+\cdots
+\theta_{n-1})
\end{array}
\right)
\end{align*}
has Jacobian determinant%
\[
r_{1}^{2n-1}t_{2}^{2n-3}t_{3}^{2n-5}\cdots t_{n-1}^{3}t_{n}.
\]
We can reduce the dimensionality of the integral for $c_{n}$ by two, using the
fact that%
\[%
{\displaystyle\int\limits_{0}^{2\pi}}
{\displaystyle\int\limits_{0}^{\infty}}
r_{1}^{2n-1}\exp\left[  -r_{1}^{2}V\left(  t_{2},t_{3},\ldots,t_{n};\theta
_{1},\theta_{2},\ldots,\theta_{n-1}\right)  \right]  dr_{1}d\theta_{0}%
=\pi(n-1)!V^{-n}.
\]
To go forward, it must be understand how $t_{i+1}$ and $\theta_{i}$ interact
with each other, as a consequence of $\max\{r_{i},r_{i+1}\}<r_{i,i+1}$.

Let $n=2$ for simplicity's sake. \ Figure 3 makes clear why $\pi/3\leq
\theta_{1}\leq5\pi/3$ is a necessary condition for $\max\{1,t_{2}^{2}%
\}\leq\left(  t_{2}\cos(\theta_{1})-1\right)  ^{2}+t_{2}^{2}\sin(\theta
_{1})^{2}$. \ A\ sufficient condition for the latter inequality is%
\[
\left\{
\begin{array}
[c]{lll}%
2\cos(\theta_{1})\leq t_{2}\leq\dfrac{1}{2\cos(\theta_{1})} &  & \text{if
}\dfrac{\pi}{3}\leq\theta_{1}<\dfrac{\pi}{2}\text{ or }\dfrac{3\pi}{2}%
<\theta_{1}\leq\dfrac{5\pi}{3}\text{,}\\
0<t_{2}<\infty &  & \text{if }\dfrac{\pi}{2}\leq\theta_{1}\leq\dfrac{3\pi}{2}%
\end{array}
\right.
\]
and here is a proof. \ The right hand side reduces to $t_{2}^{2}-2t_{2}%
\cos(\theta_{1})+1$; when $t_{2}=2\cos(\theta_{1})$, this is equal to
$t_{2}^{2}-t_{2}^{2}+1=1$; when $t_{2}=1/(2\cos(\theta_{1}))$, this is equal
to $t_{2}^{2}-1+1=t_{2}^{2}$. On the one hand, at an interior point $t_{2}=1$,
the RHS is equal to $2\left(  1-\cos(\theta_{1})\right)  $ and is $\geq1$ iff
$\cos(\theta_{1})\leq1/2$, which is always true in our domain. \ On the other
hand, at a left exterior point $t_{2}=\cos(\theta_{1})>0$, the RHS\ is equal
to $1-\cos(\theta_{1})^{2}$ and is \ $<1$; \ at a right exterior point
$t_{2}=1/\cos(\theta_{1})>0$, the RHS\ is equal to $1/\cos(\theta_{1})^{2}-1$
and is \ $<t_{2}^{2}$. \ Finally, if $\cos(\theta_{1})\leq0$, then the RHS\ is
$\geq t_{2}^{2}+1$ and is easily\ $\geq\max\{1,t_{2}^{2}\}$.%
\begin{figure}[ptb]%
\centering
\includegraphics[
height=4.683in,
width=5.7666in
]%
{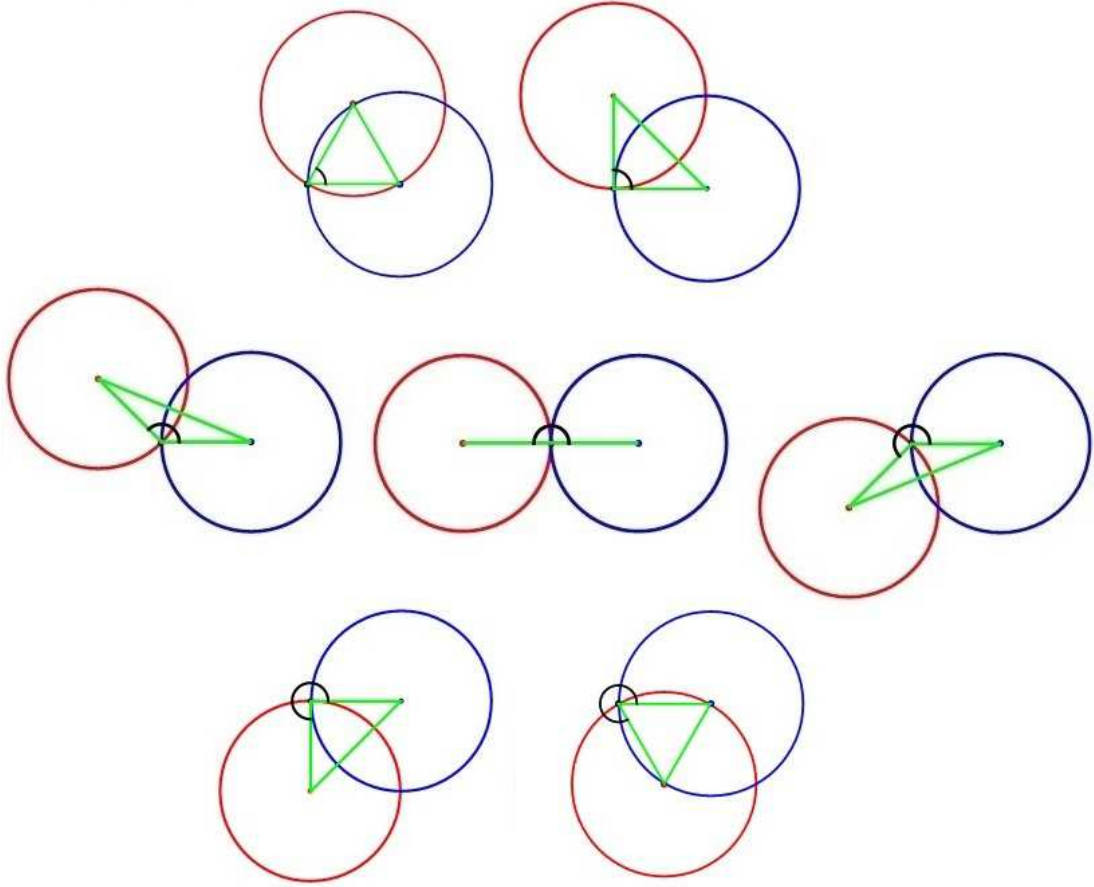}%
\caption{Top row: $\pi/3\leq\theta_{1}\leq\pi/2$. \ Middle row: $\pi
/2<\theta_{1}<3\pi/2$. \ Bottom row: $3\pi/2\leq\theta_{1}\leq5\pi/3$.
\ $t_{2}=1$ throughout.}%
\end{figure}

A\ less formal verification of sufficiency is provided in Figures 4 and 5,
illustrating the cases when $\theta_{1}$ is the midpoint of the interval
$[\pi/3,\pi/2)$ and of $(3\pi/2,5\pi/3]$ respectively. \ A trivial case is
$\theta_{1}=\pi$ since the two disks are horizontally tangent at the origin
and hence centers are automatically remote, regardless of the value of $t_{2}%
$.%
\begin{figure}[ptb]%
\centering
\includegraphics[
height=5.4656in,
width=5.866in
]%
{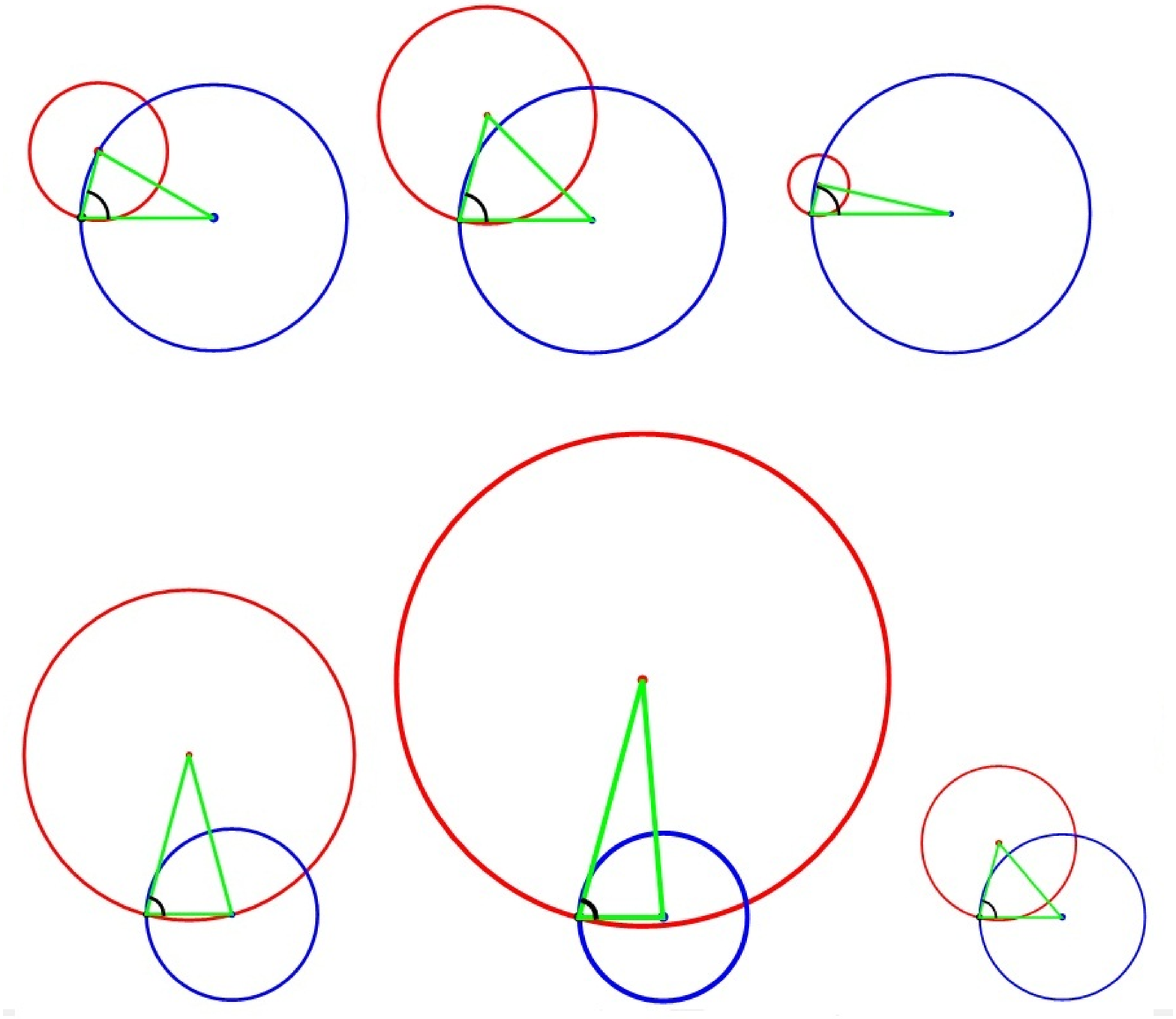}%
\caption{Top row: $2\cos(\theta_{1})$ is smallest length for which disks have
remote centers (violated in third configuration). \ Bottom row: $1/(2\cos
(\theta_{1}))$ is largest length for which disks have remote centers (violated
in second configuration). \ $\theta_{1}=5\pi/12$ throughout.}%
\end{figure}
\begin{figure}[ptb]%
\centering
\includegraphics[
height=5.2226in,
width=5.9456in
]%
{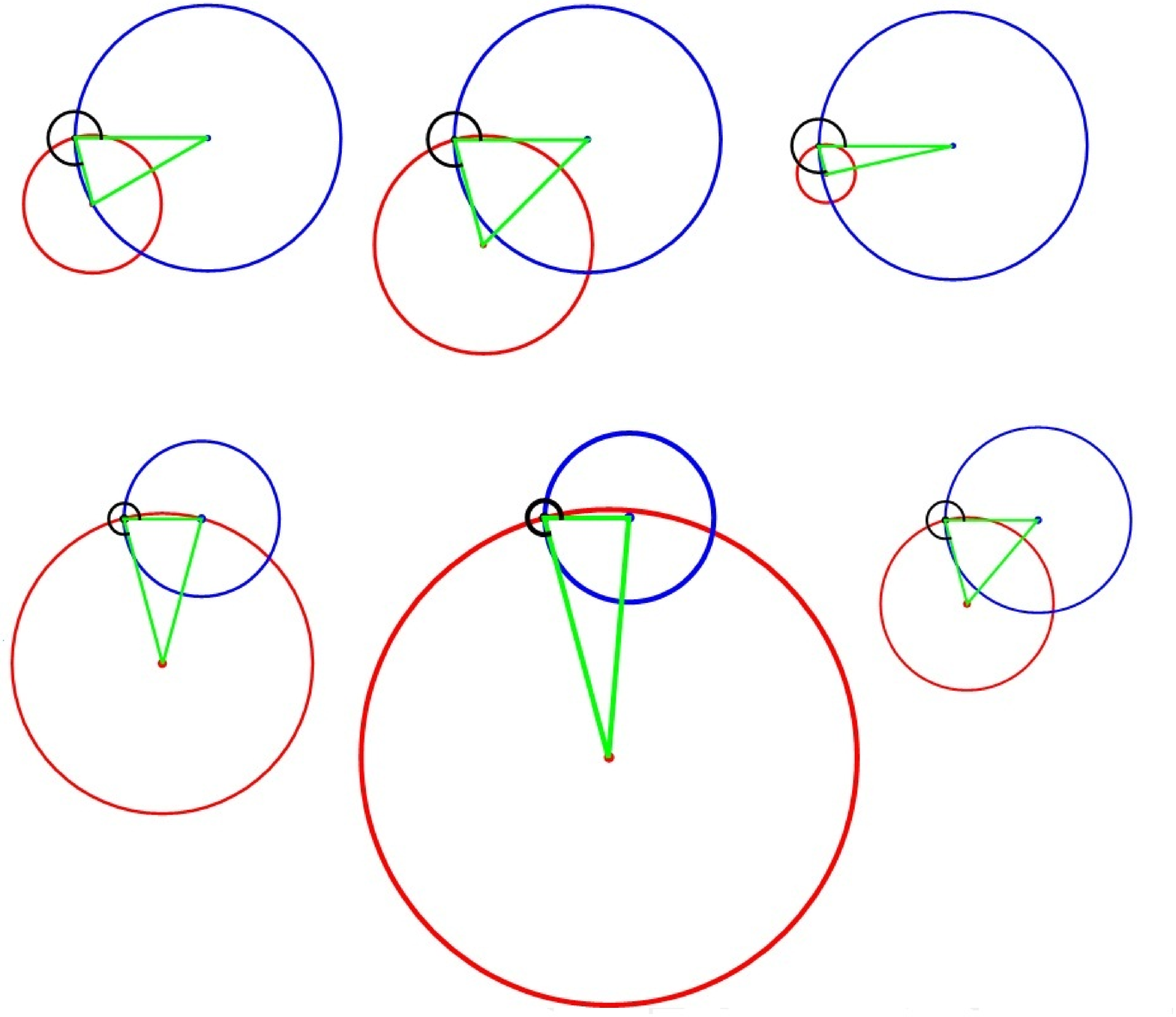}%
\caption{Top row: $2\cos(\theta_{1})$ is smallest length for which disks have
remote centers (violated in third configuration). \ Bottom row: $1/(2\cos
(\theta_{1}))$ is largest length for which disks have remote centers (violated
in second configuration). \ $\theta_{1}=19\pi/12$ throughout.}%
\end{figure}

We therefore have%
\begin{align*}
c_{2}  & =\frac{\pi}{2}\left(
{\displaystyle\int\limits_{\pi/3}^{\pi/2}}
\;%
{\displaystyle\int\limits_{2\cos(\theta_{1})}^{1/(2\cos(\theta_{1}))}}
+%
{\displaystyle\int\limits_{\pi/2}^{3\pi/2}}
\;%
{\displaystyle\int\limits_{0}^{\infty}}
+%
{\displaystyle\int\limits_{3\pi/2}^{5\pi/3}}
\;%
{\displaystyle\int\limits_{2\cos(\theta_{1})}^{1/(2\cos(\theta_{1}))}}
\right)  \frac{t_{2}}{V(t_{2};\theta_{1})^{2}}dt_{2}d\theta_{1}\\
& =0.316585....
\end{align*}
Another representation involves a binary variable $\sigma_{i}$ defined as
\[
\sigma_{i}=\left\{
\begin{array}
[c]{lll}%
0 &  & \text{if }\theta_{i}>\pi/2,\\
1 &  & \text{if }\theta_{i}<\pi/2.
\end{array}
\right.
\]
Let $I(\sigma_{1},\sigma_{2},\ldots,\sigma_{n})$ be the contribution to
$c_{n}$ with $\{\theta_{1},\theta_{2},\ldots,\theta_{n}\}$ in the range
specified by $\sigma_{1},\sigma_{2},\ldots,\sigma_{n}$. \ Clearly%
\[
c_{2}=I(1,0)+I(0,0)+I(0,1)=I(0,0)+2I(1,0)
\]
and, more generally,%
\[
c_{3}=I(0,0,0)+3I(1,0,0)+3I(1,1,0),
\]%
\[
c_{4}=4I(0,0,0,0)+4I(1,1,0,0)+2I(1,0,1,0)+4I(1,1,1,0),
\]%
\[
c_{5}=5I(1,1,1,1,0)+I(1,1,1,1,1).
\]
\{The coefficients in formulas (19)\ and\ (20) of \cite{TW-DsUn} should be
$\pi/2$, not $\pi$. \ Tao \&\ Wu's numerical result $0.316333...$ agrees with
ours to three decimal places.\} \ From%
\[
I(1,1,0)=\frac{2\pi}{3}%
{\displaystyle\int\limits_{\pi/3}^{\pi/2}}
\;%
{\displaystyle\int\limits_{\pi/3}^{\pi/2}}
\;%
{\displaystyle\int\limits_{2\cos(\theta_{1})}^{1/(2\cos(\theta_{1}))}}
\;%
{\displaystyle\int\limits_{2\cos(\theta_{2})}^{1/(2\cos(\theta_{2}))}}
\frac{t_{2}^{3}t_{3}}{V(t_{2},t_{3};\theta_{1},\theta_{2})^{3}}dt_{3}%
dt_{2}d\theta_{2}d\theta_{1},
\]%
\[
I(1,0,0)=\frac{2\pi}{3}%
{\displaystyle\int\limits_{\pi/3}^{\pi/2}}
\;%
{\displaystyle\int\limits_{\pi/2}^{3\pi/2-\theta_{1}}}
\;%
{\displaystyle\int\limits_{2\cos(\theta_{1})}^{1/(2\cos(\theta_{1}))}}
\;%
{\displaystyle\int\limits_{0}^{\infty}}
\frac{t_{2}^{3}t_{3}}{V(t_{2},t_{3};\theta_{1},\theta_{2})^{3}}dt_{3}%
dt_{2}d\theta_{2}d\theta_{1},
\]%
\[
I(0,0,0)=\frac{2\pi}{3}%
{\displaystyle\int\limits_{\pi/2}^{\pi}}
\;%
{\displaystyle\int\limits_{\pi/2}^{3\pi/2-\theta_{1}}}
\;%
{\displaystyle\int\limits_{0}^{\infty}}
\,%
{\displaystyle\int\limits_{0}^{\infty}}
\frac{t_{2}^{3}t_{3}}{V(t_{2},t_{3};\theta_{1},\theta_{2})^{3}}dt_{3}%
dt_{2}d\theta_{2}d\theta_{1}%
\]
we obtain $c_{3}=0.033056...$. \ Note that $\theta_{3}>\pi/2$ implies that
$\theta_{2}=2\pi-\theta_{3}-\theta_{1}<3\pi/2-\theta_{1}$, as indicated for
both $I(1,0,0)$ and $I(0,0,0)$. \ \ \{Tao \&\ Wu's numerical result
$0.032939...$ agrees with ours to two decimal places.\} \ What is missing,
however, is proof that some hitherto undetected interaction between $t_{2}$,
$t_{3}$, $\theta_{1}$, $\theta_{2}$ does not exist. \ Figure 6 exhibits
underlying complexity; we regrettably must stop here.%
\begin{figure}[ptb]%
\centering
\includegraphics[
height=6.4792in,
width=6.0217in
]%
{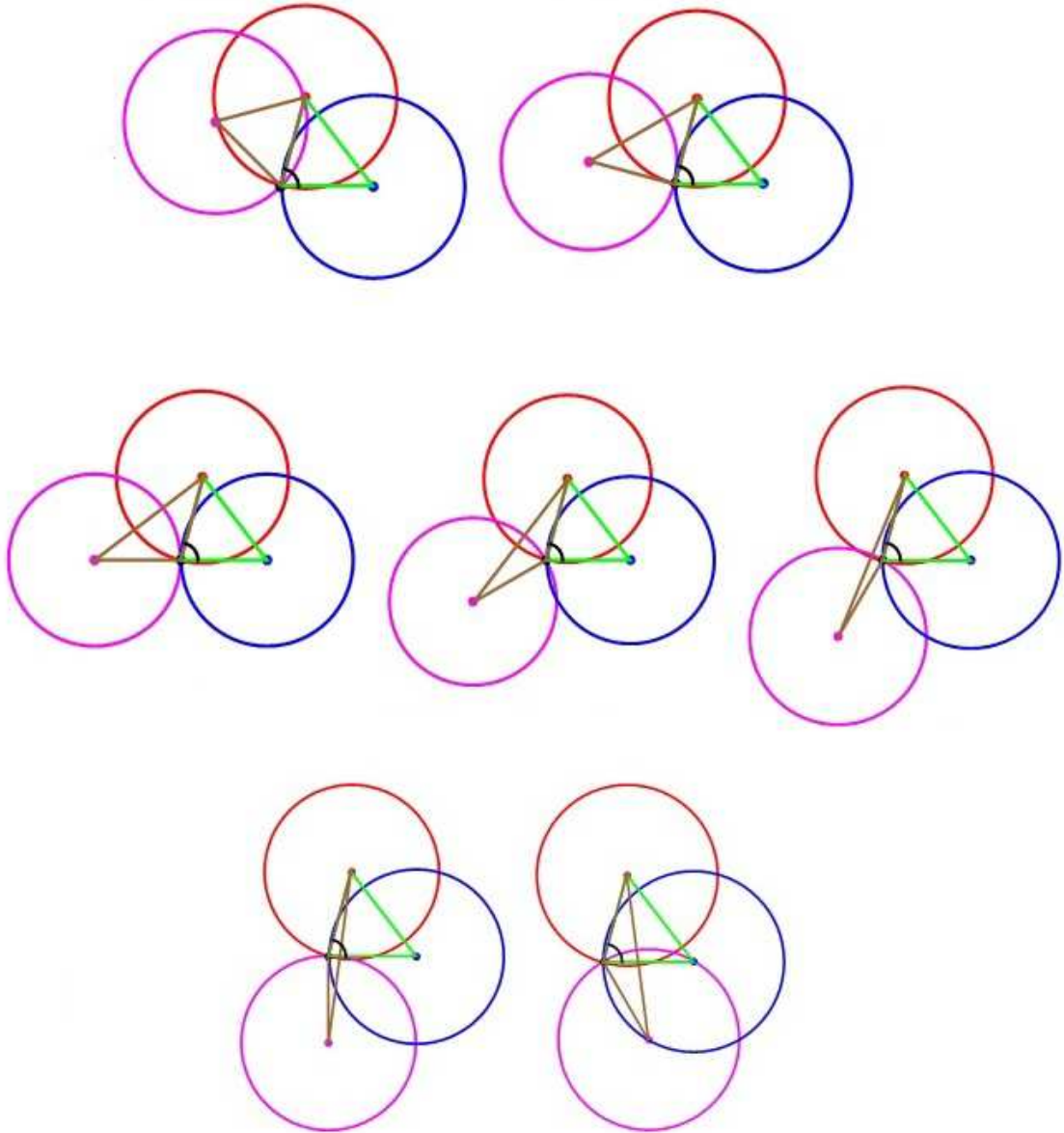}%
\caption{Top row: $\pi/3\leq\theta_{2}\leq\pi/2$. \ Middle row: $\pi
/2<\theta_{2}<13\pi/12$. \ Bottom row: $13\pi/12\leq\theta_{2}\leq5\pi/4$.
\ These correspond to $I(1,1,0)$, $I(1,0,0)$, $I(1,0,1)$ respectively.
\ $t_{2}=t_{3}=1$ and $\theta_{1}=5\pi/12$ throughout.}%
\end{figure}

\section{Acknowledgements}

I\ am grateful to Emeritus Professor Fa-Yueh Wu (Northeastern) and Kevin Daily
(Wolfram Research) for kind and helpful correspondence. \ On some future day,
if the \textsc{Mathematica 10} commands
\[%
\begin{array}
[c]{l}%
\text{\texttt{NIntegrate[Boole[}}\\
\text{\texttt{x1\symbol{94}2+y1\symbol{94}2
$<$%
= (x1-x2)\symbol{94}2+(y1-y2)\symbol{94}2 \&\& }}\\
\text{\texttt{x2\symbol{94}2+y2\symbol{94}2
$<$%
= (x1-x2)\symbol{94}2+(y1-y2)\symbol{94}2] *}}\\
\text{\texttt{Exp[-RegionMeasure[RegionUnion[}}\\
\text{\texttt{Disk[\{x1,y1\}, Sqrt[x1\symbol{94}2+y1\symbol{94}2]],}}\\
\text{\texttt{Disk[\{x2,y2\}, Sqrt[x2\symbol{94}2+y2\symbol{94}2]]]]],}}\\
\text{\texttt{\{x1,-Infinity,Infinity\},\{y1,-Infinity,Infinity\},}}\\
\text{\texttt{\{x2,-Infinity,Infinity\},\{y2,-Infinity,Infinity\},}}\\
\text{\texttt{Exclusions -%
$>$
\{x1\symbol{94}2+y1\symbol{94}2 == 0
$\vert$%
$\vert$%
}}\\
\text{\texttt{x2\symbol{94}2+y2\symbol{94}2 == 0\}]}}%
\end{array}
\]

\noindent and%
\[%
\begin{array}
[c]{l}%
\text{\texttt{NIntegrate[Boole[}}\\
\text{\texttt{x1\symbol{94}2+y1\symbol{94}2
$<$%
= (x1-x2)\symbol{94}2+(y1-y2)\symbol{94}2 \&\& }}\\
\text{\texttt{x2\symbol{94}2+y2\symbol{94}2
$<$%
= (x1-x2)\symbol{94}2+(y1-y2)\symbol{94}2 \&\& }}\\
\text{\texttt{x1\symbol{94}2+y1\symbol{94}2
$<$%
= (x1-x3)\symbol{94}2+(y1-y3)\symbol{94}2 \&\& }}\\
\text{\texttt{x3\symbol{94}2+y3\symbol{94}2
$<$%
= (x1-x3)\symbol{94}2+(y1-y3)\symbol{94}2 \&\& }}\\
\text{\texttt{x2\symbol{94}2+y2\symbol{94}2
$<$%
= (x2-x3)\symbol{94}2+(y2-y3)\symbol{94}2 \&\& }}\\
\text{\texttt{x3\symbol{94}2+y3\symbol{94}2
$<$%
= (x2-x3)\symbol{94}2+(y2-y3)\symbol{94}2] *}}\\
\text{\texttt{Exp[-RegionMeasure[RegionUnion[}}\\
\text{\texttt{Disk[\{x1,y1\}, Sqrt[x1\symbol{94}2+y1\symbol{94}2]],}}\\
\text{\texttt{Disk[\{x2,y2\}, Sqrt[x2\symbol{94}2+y2\symbol{94}2]],}}\\
\text{\texttt{Disk[\{x3,y3\}, Sqrt[x3\symbol{94}2+y3\symbol{94}2]]]]],}}\\
\text{\texttt{\{x1,-Infinity,Infinity\},\{y1,-Infinity,Infinity\},}}\\
\text{\texttt{\{x2,-Infinity,Infinity\},\{y2,-Infinity,Infinity\},}}\\
\text{\texttt{\{x3,-Infinity,Infinity\},\{y3,-Infinity,Infinity\},}}\\
\text{\texttt{Exclusions -%
$>$
\{x1\symbol{94}2+y1\symbol{94}2 == 0
$\vert$%
$\vert$%
}}\\
\text{\texttt{x2\symbol{94}2+y2\symbol{94}2 == 0
$\vert$%
$\vert$%
}}\\
\text{\texttt{x3\symbol{94}2+y3\symbol{94}2 == 0\}]}}%
\end{array}
\]

\noindent become numerically feasible, then Section 2 of this paper will be
rendered unnecessary.

\end{document}